\documentclass{amsproc}%
\usepackage{amsfonts}
\usepackage{amsmath}
\usepackage{amssymb}
\usepackage{graphicx}%
\providecommand{\U}[1]{\protect\rule{.1in}{.1in}}

\theoremstyle{plain}

\newtheorem{proposition}{Proposition}

\newtheorem{theorem}{Theorem}
\numberwithin{equation}{section}
\begin{document}
\title[Wave equation]{The reconstruction of a wave equation from one side measurement}
\author{Amin Boumenir}
\author{Vu Kim Tuan}
\address{Department of Mathematics,\\
University of West Georgia, GA 30118, USA}
\email{boumenir@westga.edu, vu@westga.edu}
\subjclass[2000]{Primary 35R30, 35L05; Secondary 34A55 }
\keywords{Wave equation, boundary inversion, Gelfand-Levitan theory}

\begin{abstract}
We are concerned with the reconstruction of a one dimensional wave equation,
where the potential is known in a neighborhood of one of the end points of the boundary. We show then the sought potential can be determined by one single measurement of the solution at that end. We also show that the solution can be evaluated at the other end without the need of solving the wave equation.
\end{abstract}
\maketitle

\section{Introduction}

Consider the problem of recovering an unknown coefficient $q$ and constants
$\ell, H$, appearing in the following wave equation
\begin{equation}
\left\{
\begin{array}
[c]{l}%
u_{tt}(x,t)=u_{xx}(x,t)-q(x)u(x,t),\ \ \ \ \ \ \ 0<x<\ell,\text{\ \ }\ t>0,\\
u_{x}\left(  0,t\right)  -hu(0,t)=0,\text{ \ \ \ \ }u_{x}\left(
\ell,t\right)  +Hu(\ell,t)=0,\quad t>0,\\
u(x,0)=f(x)\ \ \ \ \ \text{and} \ \ \ \ \ u_{t}\left(  x ,0\right)
=g(x),\quad0<x< \ell.
\end{array}
\right.\label{wa}
\end{equation}
This note is concerned with choosing initial conditions $\{ f,g \}$ such that
the data gathered from a single measurement of $u(0,t)$ for $t>0$ is enough
for the reconstruction of the real valued function $q\in L\left(  0,\ell\right)
$. This will answer the question of how to recover $q$ when one of the endpoints
on the boundary is not accessible for measurements, for example when dealing
with graphs \cite{advwav, boutuanwav}, or clogged buried pipelines or blood
vessels. To this end we require that%
\begin{equation}
h\in\mathbb{R}\;\text{ is given, and}\;\;q(x)\ \ \text{is known for}\;0<
x<\epsilon,\;\;\;\mathrm{where}\;\;\;\epsilon
\;\;\mathrm{is\;arbitrarily\;small.} \label{q}%
\end{equation}
In \cite{boutuanwav}, it was shown that at most two measurements on both sides
were required, i.e. two pairs of traces $\left\{  u(0,t),u(\ell,t)\right\}  $
that were generated by two special initial conditions. The issue here is to
choose initial conditions $\left\{  f,\ g\right\}  $ in (\ref{wa}) such that
the observation $u(0, t)$ carries all the necessary spectral data in order to
apply the Gelfand-Levitan theory.

\textbf{Statement of the problem:} Assuming that (\ref{q}) holds, can we
reconstruct $q$, $\ell,$ and $H$ from one single measurement of $\left\{
u\left(0, t\right)  \right\}_{t>0} $? Also can we compute  $u(\ell,t)$  using  $u(0,t)$ and $\{f,g\}$ only?

The main motivation behind the above setting is the problem of unclogging an
inaccessible pipe such as a buried pipeline or a blood vessel, that can be
accessed from one end only. The idea here is to locate the obstruction, i.e. $\ell
$, and its type $H$, and then generate a pressure wave $u(\ell,t)$ that can dissolve it. This will avoid constructing $q$ and then solving the wave equation (\ref{wa}).

\section{Preliminaries}

The solution of (\ref{wa}) can be written as
\begin{equation}
u(x,t)=\sum_{n\geq0}\left(  a(\lambda_{n})\cos(t\sqrt{\lambda_{n}}%
)+b(\lambda_{n})\frac{\sin(t\sqrt{\lambda_{n}})}{\sqrt{\lambda_{n}}}\right)
\frac{1}{\alpha_{n}^{2}}y(x,\lambda_{n}), \label{w}%
\end{equation}
where the eigenfunctions $y(x,\mathbb{\lambda}_{n})$, for $n=0,1,2,\cdots,$ are solutions of
\[
\left\{
\begin{array}
[c]{l}%
-y^{\prime\prime}(x,\lambda)+q(x)y(x,\lambda)=\lambda y(x,\lambda),\quad0<x< \ell,\\
y(0,\lambda)=1,\;\;y^{\prime}(0,\lambda)=h,
\end{array}
\right.
\] and the eigenvalues, $\mathbb{\lambda}_{0},\mathbb{\lambda}_{1},\cdots,$ are
determined as the roots of
$ y^{\prime}(\ell,\lambda)+Hy(\ell,\lambda)=0$. Denote  their norms by $ \quad\alpha_{n}=\left\Vert y\left(x,\lambda_{n}\right)\right\Vert_{L^{2}(0,\ell)}$ and their Fourier coefficients  by
\[
a\left(  \lambda_{n}\right)  =\int_{0}^{\ell}f(x) y(x,\lambda_{n})dx, \quad
b\left(  \lambda_{n}\right)  =\int_{0}^{\ell} g(x) y(x,\lambda_{n})dx.\]
If we assume that we are given $h$ and $q(x)$ for $0<x<\epsilon$, then we
would know $y(x,\mathbb{\lambda})$ over the interval $[0,\epsilon]$. The
transmutation operators, \cite{gas}, would also be known over the same
interval
\begin{equation}
\cos(x\sqrt{\mathbb{\lambda}}) =y(x,\mathbb{\lambda})+\int_{0}^{x}H(x,t)y(x,\mathbb{\lambda})dt,\;\;\;\;\;\;\;\quad0<x<\epsilon. \label{tra}%
\end{equation}
In order to extract the complete spectral data $\left\{\lambda_n,\alpha_n^2\right\}$ from (\ref{w}) we need to ensure all the Fourier coefficients to satisfy $a^2(\lambda_n)+b^2(\lambda_n)\neq 0$. For simplicity we show  the following proposition
\begin{proposition}
Given $q$ over $(0,\epsilon)$, where $0<\epsilon<1$ and $h\in\mathbb{R},$ we can construct a continuous function $g_{\epsilon}$ explicitly, such that\\[-3mm]\[
supp(g_{\epsilon})\subset\lbrack0,\epsilon]\;\;\;\mathrm{and}\;\;\;\int
_{0}^{\ell}g_{\epsilon}(x)y(x,\mathbb{\lambda})dx\neq 0\;\;\;\mathrm{for\;\;all}\;\;\;\mathbb{\lambda}\in\mathbb{R},
\] and in particular we have $b(\lambda_n)\neq 0$ for all  $n\geq 0$.
\end{proposition}
Proof. Let $\psi(x)=(x-\epsilon)^21_{[0,\epsilon]}(x)$, where
$1_{[0,\epsilon]}(x)=1$ if $x\in\lbrack0,\epsilon]$ and $0$ otherwise. Its
cosine transform, given by
\begin{equation}
\int_{0}^{\ell}\psi(x)\cos(x\sqrt{\lambda})dx=\int_{0}^{\epsilon}%
(\epsilon-x)^2\cos(x\sqrt{\lambda})dx=2\frac{\epsilon\sqrt{\lambda}-\sin(\epsilon\sqrt{\lambda})}{\lambda\sqrt{\lambda}}=O\left(
\frac{1}{\lambda}\right)  ,\quad\lambda\rightarrow \infty, \label{f1}%
\end{equation}
does not vanish for any $\lambda\in\mathbb{R},$ which is easy to see, since for  $\lambda > 0$ we have 
$\epsilon\sqrt{\lambda}>\sin(\epsilon\sqrt{\lambda})$, and  if $\lambda \leq 0$, then the integrand is positive. Since the kernel $H(x,t)$, in (\ref{tra}), is known for $0<t\leq x<\epsilon$, we also have
\begin{align*}
0\neq\int_{0}^{\ell}\psi(x)\cos(x\sqrt{\lambda})dx  & 
 =\int_{0}^{\epsilon}\psi(x)y(x,\lambda)dx+\int_{0}^{\epsilon}
 \int_{0}^{x}H(x,t)y(t,\mathbb{\lambda})dt\psi(x)dx\\
&  =\int_{0}^{\epsilon}\left\{  \psi(t)+\int_{t}^{\epsilon}H(x,t)\psi
(x)dx\right\}  y(t,\mathbb{\lambda})dt.
\end{align*}
Thus we can define the initial condition
\begin{equation}
g_{\epsilon}(x)=\left\{
\begin{array}
[c]{ll}%
\psi(x)+\int_{x}^{\epsilon}H(t,x)\psi(t)dt, & 0<x<\epsilon\\
0, & \epsilon\leq x<\ell
\end{array}
\right.  , \label{ge}%
\end{equation}
which vanishes for $x>\epsilon$ and its $y$-transform, by (\ref{f1}), has no
zeros for $\mathbb{\lambda}\in\mathbb{R}$. This initial condition will
essentially ensure that the Fourier coefficients $b(\lambda_n)$ do not vanish and are known
explicitly
\begin{equation}
b(\lambda_{n})=\int_{0}^{\ell}g_{\epsilon}(x)y(x,\mathbb{\lambda}_{n})dx
=2\frac{\epsilon\sqrt{\lambda_n}-\sin(\epsilon\sqrt{\lambda_n})}{\lambda_n\sqrt{\lambda_n}}\neq
0,\;\;\mathrm{and}\;\;b(\mathbb{\lambda}_{n})=O\left(  \dfrac{1}%
{\mathbb{\lambda}_{n}}\right)  . \label{bn}%
\end{equation}
As for the initial condition $f$ in (\ref{wa}), we take $f(x)=0$, and so
$a(\lambda_{n})=0$.

\section{Reading the data at $x=0$}

We now use the asymptotics of the eigenvalues $\lambda_{n} \sim\left(
\frac{\pi}{\ell} n\right)^{2}$, \cite{mar}, of the norming constants 
$\alpha_{n}^{2}\sim\frac{\ell}{2}$, \cite{mar}, and of the Fourier coefficients 
$b(\lambda_{n}) = O\left(\frac{1}{\lambda_{n}}\right) = O\left(\frac{1}{n^{2}}\right)$ of
the initial condition defined by (\ref{ge}), to deduce that the observation
\begin{equation}
u(0,t)=2\sum_{n\geq0}\frac{\epsilon\sqrt{\lambda_n}-\sin(\epsilon\sqrt{\lambda_n})}{\lambda^2_n}\frac{\sin(t\sqrt{\lambda_{n}})}{\alpha_{n}^{2}}\ \label{w1}%
\end{equation}
is a continuous function in $t$, and by the Lebesgue dominated convergence
theorem one can apply the Laplace transform to \eqref{w1} termwise to obtain
\begin{equation}
\mathcal{L}(u)(0,s)=\sum_{n\geq0}\frac{\epsilon\sqrt{\lambda_n}-\sin(\epsilon\sqrt{\lambda_n})}{\lambda^2_n}\frac{2\sqrt{\lambda_n}}{\left(  s^{2}%
+\mathbb{\lambda}_{n}\right) \alpha_{n}^{2}}, \quad\Re(s) > 0. \label{w2}%
\end{equation}
The series in \eqref{w2} represents a meromorphic function, and one can use
pole finding methods to read off all $\left\{  \pm i\sqrt{\mathbb{\lambda}%
_{n}}\right\}  $, i.e. all the eigenvalues that appear in (\ref{w1}) and then compute the
residues of $\mathcal{L}(u)(0,s)$ at $\left\{ i\sqrt{\mathbb{\lambda}%
_{n}}\right\}  $ to obtain the sequence
\begin{equation}
\left\{ \frac{\epsilon\sqrt{\lambda_n}-\sin(\epsilon\sqrt{\lambda_n})}{i\lambda^2_n\alpha_{n}^{2}}\right\}_{n\geq0}. \label{bb}%
\end{equation}
Recall that  $\left\{
\lambda_{n}\right\}$ are known from  (\ref{w2}),
and so we can deduce all the norming constants $\left\{
\alpha_{n}^{2}\right\}_{n\geq0}$ from (\ref{bb}). In other words we have the complete
spectral data $\left\{  \lambda_{n}, \alpha_{n}^{2}\right\}  $ that is
required to form the spectral function and so $q$, $\ell$ and $H$ by the Gelfand-Levitan theory, \cite{gas}. It is known
that asymptotics of the eigenvalues have the form \cite{mar}
\begin{equation}
\sqrt{\lambda_{n}}=\frac{\pi}{\ell}n+\dfrac{1}{n\ell}a_{1}+o\left(  \dfrac
{1}{n}\right),\;\;\text{where}\;\; a_{1}=  h+H+\frac{1}{2}\int_{0}^{\ell
}q(x)dx.\end{equation} Thus taking the limits, we find
\[
\ell=\lim_{n\rightarrow\infty}\dfrac{n\pi}{\sqrt{\lambda_{n}}}%
\;\;\;\;\;\;\mathrm{and}\;\;\;\;\;\; \lim_{n\rightarrow\infty} n\left(
\ell\sqrt{\mathbb{\lambda}_{n}}-n\pi\right)  =h+H+\frac{1}{2}\int_{0}^{\ell
}q(x)dx,
\]
  and as we now know $q,\;\ell$ and $h$, we can get $H$. Thus we
have proved

\begin{theorem}
Assume that we know $h\in\mathbb{R}$ and $q(x)$ for $0<x<\epsilon$, with
$\epsilon$ arbitrarily small. Then using the initial conditions $f=0$,
$g=g_{\epsilon}$ as in \eqref{ge}, we can uniquely reconstruct $q$ over
$(0,\ell)$ and $H$, from one single reading of $u(0,t)$ for $t>0$.
\end{theorem}

We next examine the limiting case $h\rightarrow\infty$, which is the Dirichlet
boundary condition at $x=0.$

\section{Dirichlet condition}

We now briefly show that the above approach extends to the Dirichlet case
\[
\left\{
\begin{array}
[c]{l}%
u_{tt}(x,t)=u_{xx}(x,t)-q(x)u(x,t),\ \ \ \ \ \ \ 0<x<\ell,\text{\ \ }\ t>0,\\
u(0,t)=0,\text{ \ \ \ \ \ \ \ \ \ \ \ \ }u_{x}\left(  \ell,t\right)
+Hu(\ell,t)=0,\quad t>0,\\
u(x,0)=f(x)\ \ \ \ \ \text{and}\ \ \ \ \ u_{t}\left(  x,0\right)
=g(x),\quad0<x<\ell.
\end{array}
\right.
\]
Since $u(0,t) = 0$, we measure $u_x(0,t)$ for $t>0$, which is given by
\begin{equation}
u_x(x,t)=\sum_{n\geq0}\left(  \widetilde{a}(\mu_{n})\cos(t\sqrt{\mu_{n}%
})+\widetilde{b}(\mu_{n})\frac{\sin(t\sqrt{\mu_{n}})}{\sqrt{\mu_{n}}}\right)
\frac{1}{\widetilde{\alpha}_{n}^{2}}\varphi^{\prime}(x,\mu_{n}),
\end{equation}
where the eigenfunctions $\varphi(x,\mathbb{\mu}_{n})$, for $n=0,1,2,\cdots,$
are solutions of
\[
\left\{
\begin{array}
[c]{l}%
-\varphi^{\prime\prime}(x,\mu)+q(x)\varphi(x,\mu)=\mu\varphi(x,\mu),\quad0<x<\ell,\\
\varphi(0,\mu)=0\;\;\text{and}\;\;\varphi^{\prime}(0,\mu)=1.
\end{array}
\right.
\] The eigenvalues $\mu_{n}$ are the roots of
\begin{equation}
\varphi^{\prime}(\ell,\mu_{n})+H\varphi(\ell,\mu_{n})=0,
\end{equation} and the Fourier coefficients and norming constants are 
\[
\widetilde{a}\left(  \mu_{n}\right)  =\int_{0}^{\ell}  f(x) \varphi(x,\mu_{n})dx,\quad\widetilde{b}\left( \mu_{n}\right)  =\int_{0}^{\ell
}g(x) \varphi(x,\mu_{n})dx,\quad\widetilde{\alpha}_{n}=\left\Vert
\varphi\left(  x,\mu_{n}\right)  \right\Vert _{L^{2}(0,\ell)}.
\]
If we again assume that we are given $q(x)$ for $0<x<\epsilon$, then we would
know $\varphi(x,\mathbb{\mu})$ over the interval $[0,\epsilon]$ and so the
transmutation operators, \cite{levbook}, is known explicitly
\begin{equation}
\frac{\sin(x\sqrt{\mathbb{\mu}})}{\sqrt{\mu}}=\varphi(x,\mathbb{\mu})+\int_{0}^{x}L(x,t)\varphi(t,\mathbb{\mu})dt,\;\;\;\;\;\;\quad0<x<\epsilon.\label{tra2}
\end{equation}
We now need the following proposition

\begin{proposition}
Given $q$ over $(0,\epsilon)$, where $0<\epsilon<1,$ then we can construct a
continuous function $\widetilde{g}_{\epsilon}$ explicitly, such that
\[
supp(\widetilde{g}_{\epsilon})\subset\lbrack0,\epsilon]\;\;\;\mathrm{and}%
\;\;\;\int_{0}^{\ell}\widetilde{g}_{\epsilon}(x)\varphi(x,\mathbb{\mu}%
)dx\neq0\;\;\;\mathrm{for\;\;all}\;\;\;\mathbb{\mu}\in\mathbb{R}.
\]
\end{proposition}

Proof. If we use again ${\psi}(x)=(\epsilon - x)1_{[0,\epsilon]}(x)$, then its
sine transform, given by
\[
\int_{0}^{\ell}{\psi}(x)\frac{1}{\sqrt{\mu}}\sin(x\sqrt{\mu})dx=\int
_{0}^{\epsilon}(\epsilon-x)\frac{1}{\sqrt{\mu}}\sin(x\sqrt{\mu})dx=\frac{\epsilon \sqrt{\mu} - \sin \left( \epsilon \sqrt{\mu}\right)}{\mu}=O\left(\frac{1}{\sqrt{\mu}}\right)  ,\quad\mu \rightarrow \infty,
\]
does not vanish for any $\mu\in\mathbb{R}$. Since the kernel $L(x,t)$
in (\ref{tra2}) is known for $0<t\leq x<\epsilon$, we also have
\begin{align}
0\neq\int_{0}^{\ell}{\psi}(x)\frac{1}{\sqrt{\mu}}\sin(x\sqrt{\mu})dx  &
=\int_{0}^{\epsilon}{\psi}(x)\varphi(x,\mu)dx+\int_{0}^{\epsilon}\int_{0}%
^{x}L(x,t)\varphi(t,\mathbb{\mu})dt {\psi}(x)dx\nonumber\\
&  =\int_{0}^{\epsilon}\left\{ { \psi}(t)+\int_{t}^{\epsilon}L(x,t){\psi}
(x)dx\right\}  \varphi(t,\mathbb{\mu})dt.\label{f12}
\end{align}
Thus, we define the initial condition
\begin{equation}
\widetilde{g}_{\epsilon}(x)=\left\{
\begin{array}
[c]{ll}%
{\psi}(x)+\int_{x}^{\epsilon}L(t,x){\psi}(t)dt, & 0<x<\epsilon,\\
0, & \epsilon\leq x<\ell,
\end{array}
\right.  \label{gd}
\end{equation}
which vanishes for $x>\epsilon$ and its $\varphi$-transform, by (\ref{f12}),
has no zeros for $\mathbb{\mu}\in\mathbb{R}$. This initial condition will
essentially ensure that the Fourier coefficients $\widetilde{b}(\mu_{n})$ do not vanish and are known
explicitly
\[
\widetilde{b}(\mu_{n})=\int_{0}^{\ell}\widetilde{g}_{\epsilon}(x)\varphi
(x,\mathbb{\mu}_{n})dx=\frac{\epsilon \sqrt{\mu_{n}}- \sin(\epsilon\sqrt{\mu_{n}})}{\mu_{n}} \neq0.
\]
Taking $f = 0$, the observation $u_x(0,t)$ has the form
\begin{equation}
u_x(0,t)=\sum_{n\geq0}\frac{\epsilon \sqrt{\mu_n}- \sin(\epsilon\sqrt{\mu_n})}{{\mu_n}^{3/2}} \frac{\sin(t\sqrt{\mu_{n}})}{\widetilde{\alpha}_{n}^{2}},
\end{equation}
out of which  we can extract the full spectral data $\{\mu_n, \widetilde{\alpha}_{n}\}_{n \geq 0}$.
We have just proved 
\begin{theorem}
Assume that we know $q(x)$ for $0<x<\epsilon$, with
$\epsilon$ arbitrarily small. Then using the initial conditions $f=0$,
$g=\widetilde{g}_{\epsilon}$ as in \eqref{gd}, we can uniquely reconstruct $\ell$, $q$ over
$(0,\ell)$ and $H$, from one single reading of $u_x(0,t)$ for $t>0$.
\end{theorem}
Proof. It remains to see that asymptotics for the eigenvalues $\mu_n$ are now given  by \cite{mar}
\begin{equation}
\sqrt{\mu_{n}}=\frac{\pi}{\ell}n+\dfrac{\pi}{2\ell}+\dfrac{1}{n\ell}a_{1}+o\left(  \dfrac
{1}{n}\right),\;\;\text{where}\;\; a_{1}= H+\frac{1}{2}\int_{0}^{\ell}q(x)dx,\end{equation} which means
\[
\ell=\lim_{n\rightarrow\infty}\dfrac{n\pi}{\sqrt{\mu_{n}}}%
\;\;\;\;\;\;\mathrm{and}\;\;\;\;\;\;H= \lim_{n\rightarrow\infty} n\left(
\ell\sqrt{\mathbb{\mu}_{n}}-\left(n+\dfrac{1}{2}\right)\pi\right) -\frac{1}{2}\int_{0}^{\ell
}q(x)dx.
\]

\section{Determination of $u(\ell, t)$}

We consider the case of the wave equation (\ref{wa}), where  $y(0,\lambda)=1$. Recall that recovered eigenvalues $\mathbb{\lambda}_{n}$ are the zeros of the
boundary function
\[
\Phi(\mathbb{\lambda})=y^{\prime}(\ell,\lambda)+Hy(\ell,\lambda),
\]
which is an entire function of order $1/2$, and so can also be represented by
its zeros by Hadamard factorization theorem \cite{mar}
\[
\displaystyle\Phi(\mathbb{\lambda})={\ell}(\mathbb{\lambda}_0-\mathbb{\lambda
})\prod_{n\geq1}\left(  \frac{\ell^2}{\pi^2 n^{2}}(\mathbb{\lambda}_{n}-\mathbb{\lambda})\right)  .
\]
So $\Phi(\mathbb{\lambda})$ is determined by the $\{\lambda_{n}\}_{n\geq0}$ and to
reconstruct $u(\ell,t)$ without the knowledge of $q$ we use the fact that
\[
y^{2}(x,\lambda)=\left(  y^{\prime}(x,\lambda)\partial_{\mathbb{\lambda}%
}y(x,\lambda)-y(x,\lambda)\partial_{\mathbb{\lambda}}y^{\prime}(x,\lambda
)\right)  ^{\prime},
\]
where $\partial_{\mathbb{\lambda}}$ is the derivative with respect to
$\mathbb{\lambda}$, which yields
\begin{align*}
\alpha_{n}^{2}  &  =\int_{0}^{\ell}|y(x,\lambda_{n})|^{2}dx\\
&  =y^{\prime}(\ell,\lambda_{n})\partial_{\mathbb{\lambda}}y(\ell,\lambda
_{n})-y(\ell,\lambda_{n})\partial_{\mathbb{\lambda}}y^{\prime}(\ell
,\lambda_{n})-\left(  y^{\prime}(0,\lambda_{n})\partial_{\mathbb{\lambda}%
}y(0,\lambda_{n})-y(0,\lambda_{n})\partial_{\mathbb{\lambda}}y^{\prime
}(0,\lambda_{n})\right) \\
&  =y^{\prime}(\ell,\lambda_{n})\partial_{\mathbb{\lambda}}y(\ell,\lambda
_{n})-y(\ell,\lambda_{n})\partial_{\mathbb{\lambda}}y^{\prime}(\ell
,\lambda_{n})\\
&  =\Phi(\lambda_{n})\partial_{\mathbb{\lambda}}y(\ell,\lambda_{n}%
)-y(\ell,\lambda_{n})\Phi^{\prime}(\lambda_{n})\\
&  =-y(\ell,\lambda_{n})\Phi^{\prime}(\lambda_{n}),
\end{align*}
where we used the fact that $y(0,\mathbb{\lambda})=1,\;y^{\prime
}(0,\mathbb{\lambda})=h$ and $\Phi(\mathbb{\lambda}_{n})=0$. Consequently,
\[
\frac{y(\ell,\lambda_{n})}{\alpha_{n}^{2}}=-\Phi^{\prime}(\lambda_{n}),
\]
and therefore, using (\ref{w}) and (\ref{bn}) we can deduce the profile $u(\ell,t)$
\begin{equation}
u(\ell,t)=2\sum_{n\geq0}\frac{\sin(\epsilon\sqrt{\lambda_n})-\epsilon\sqrt{\lambda_n}}{\lambda^2_n}\sin(t\sqrt{\lambda_{n}})\Phi^{\prime}(\lambda_{n})\ \label{w11}%
\end{equation}
 without the knowledge of $q$ or the  integration of the wave equation (\ref{wa}).

\end{document}